\documentclass[12pt]{article}
\usepackage{authblk}
\usepackage{mathtools} 
\usepackage{amsmath, amssymb} 
\usepackage{amsthm}
\usepackage{enumerate}  
\usepackage{url} 

\usepackage{tikz}
\usepackage{tkz-graph}
\usetikzlibrary{shapes}
\usetikzlibrary{arrows}
\usetikzlibrary{decorations.markings}
 
\usepackage{hyperref}
\usepackage{graphicx}
\usepackage{caption,subcaption}

\newcommand\cx{{\mathbb C}}

\newcommand\ints{{\mathbb Z}}
\newcommand\re{{\mathbb R}}

\newcommand\cA{{\mathcal A}}
\newcommand\cJ{{\mathcal J}}
\newcommand\cH{{\mathcal H}}

\DeclarePairedDelimiter\abs{\lvert}{\rvert}%
\DeclarePairedDelimiter\norm{\lVert}{\rVert}%

\makeatletter
\let\oldabs\abs
\def\abs{\@ifstar{\oldabs}{\oldabs*}}
\let\oldnorm\norm
\def\norm{\@ifstar{\oldnorm}{\oldnorm*}}
\makeatother

\newcommand\sbs{\subseteq}

\newcommand\pmat[1]{\begin{pmatrix} #1 \end{pmatrix}}
\newcommand\seq[4]{#1_{#2},#1_{#3},\ldots,#1_{#4}}

\newtheoremstyle{plainsl}%
	{\topsep}
	{\topsep}
	{\slshape} 
	{}
	{\normalfont\bfseries}
	{.}
	{ }
	{}

\swapnumbers

{\theoremstyle{plainsl}
\newtheorem{theorem}{Theorem}[section]
\newtheorem{lemma}[theorem]{Lemma}
}
{\theoremstyle{remark}
}

\renewcommand\proof{\noindent\textsl{Proof. }}
\newcommand\sqr[2]{{\vbox{\hrule height.#2pt
    \hbox{\vrule width.#2pt height#1pt \kern#1pt
        \vrule width.#2pt}\hrule height.#2pt}}}
\renewcommand\qed{%
	\ifmmode\eqno\sqr53
	\else\nolinebreak\ \hfill\sqr53\medbreak\fi}

\DeclareMathOperator{\rk}{rk}

\newcommand\one{{\bf1}}
\newcommand\zero{{\bf0}}

\usepackage{blkarray}

\newcommand*\pFqskip{8mu}
\catcode`,\active
\newcommand*\pFq{\begingroup
        \catcode`\,\active
        \def ,{\mskip\pFqskip\relax}%
        \dopFq
}
\catcode`\,12
\def\dopFq#1#2#3#4#5{%
        {}_{#1}F_{#2}\biggl(\genfrac..{0pt}{}{#3}{#4};#5\biggr)%
        \endgroup
}

\title{Perfect State Transfer on Weighted Graphs of the Johnson Scheme}

\author{Luc Vinet}
\author{Hanmeng Zhan}
\affil{Centre de Recherches Math\'ematiques
Universit\'e de Montr\'eal, P.O. Box 6128, Centre-ville Station,
Montr\'eal, Qu\'ebec, Canada, H3C 3J7}

\begin{document}
\maketitle

\begin{abstract}
We characterize perfect state transfer on real-weighted graphs of the Johnson scheme $\cJ(n,k)$.  Given $\cJ(n,k)=\{\seq{A}{1}{2}{k}\}$ and  $A(X) = w_0A_0 + \cdots + w_m A_m$, we show, using classical number theory results, that $X$ has perfect state transfer at time $\tau$ if and only if  $n=2k$, $m\ge 2^{\lfloor{\log_2(k)} \rfloor}$, and there are integers $\seq{c}{1}{2}{m}$ such that 
\begin{enumerate}[(i)]
\item $c_j$ is odd if and only if $j$ is a power of $2$, and 
\item for  $r=1,2,\cdots,m$, 
\[w_r = \frac{\pi}{\tau} \sum_{j=r}^m \frac{c_j}{\binom{2j}{j}} \binom{k-r}{j-r}.\]
\end{enumerate}
We then characterize perfect state transfer on unweighted graphs of $\cJ(n,k)$. In particular, we obtain a simple construction that generates all graphs of $\cJ(n,k)$ with perfect state transfer at time $\pi/2$.
\end{abstract}

\section{Introduction \label{sec_intro}}
Quantum spin networks that manifest perfect state transfer between two qubits are important devices in implementing quantum information tasks. There are two main objects of studies in the literature: spin chains with modulated nearest-neighbor couplings, and spin networks with uniform nearest-neighbor couplings. The first situation is represented by a weighted path. Here, one aims to determine weightings that allow perfect state transfer between the endpoint or internal vertices, and there have been experimental results, as well as analytic solutions based on orthogonal polynomials (see, for example, \cite{Albanese2004,Vinet2012}). In comparison, the second situation is represented by an unweighted graph. Without the freedom to change weights, the goal is to characterize graphs that allow perfect state transfer between some pair of vertices, and many results follow from algebraic graph theory (see, for example, \cite{Godsil2012,Coutinho2015}).

Given a distance regular graph, one can ``quotient" it down to a weighted path by the distance partition relative to a vertex. Hence, constructions of both types of spin networks may be related. For instance, the quotient of a hypercube is a Krawtchouk chain, and perfect state transfer between antipodal vertices of the hypercube is equivalent to perfect state transfer between the endpoints of the Krawtchouk chain (see 
Christandl, Datta, Ekert, and J. Landahl 
\cite{Christandl2004} and 
Albanese, Christandl, Datta, and Ekert
\cite{Albanese2004} for such networks that admit perfect state transfer). As another example, a Johnson graph on more than two vertices can be quotient down to a dual Hahn chain, but not one of those exhibiting perfect state transfer \cite{Albanese2004}. In fact, Ahmadi, Haghighi and Mokhtar \cite{Ahmadi2017} proved that the only generalized Johnson graphs with perfect state transfer are disjoint unions of edges. 


While most studies of perfect state transfer focus on nearest-neighbor couplings, in a realistic spin network, there may be residual couplings between non-neighboring qubits, whose strengths decrease with distance. As shown by Kay \cite{Kay2006}, it is possible to have perfect state transfer on some of these spin chains, where one allows couplings between the next-to-nearest neighbors. In a more recent paper \cite{Christandl2017}, Christandl, Vinet and Zhedanov gave the first analytic construction of next-to-nearest-neighbor Krawtchouk chains that admit perfect state transfer; moreover, some of these chains lift to weighted unions of a hypercube and its distance-$2$ graph. 

Motivated by the above findings are three questions.

\begin{enumerate}[(i)]
\item Let $X$ be an unweighted graph of diameter $d$. For $j=0,1,\cdots,d$, let $X_j$ denote the distance-$j$ graph of $X$, that is, the graph with the same vertex set as $X$ and an edge between any two vertices at distance $j$ in $X$. Can we characterize real-weighted unions of $\seq{X}{0}{1}{d}$ that admit perfect state transfer?
\item What if $X$ is distance regular? In this case, the weighted graphs enjoy many nice properties as $\{A(X_0),\cdots A(X_d)\}$ forms an association scheme.
\item What if $X$ is the Johnson graph? Any example found here yields a dual Hahn chain with couplings beyond nearest neighbors that admit perfect state transfer. 
\end{enumerate} 

This paper gives a complete answer to the third question. Using results from classical number theory and association schemes, we provide a characterization for weighted unions of the Johnson scheme $\cJ(n,k)$ that admit perfect state transfer. As a consequence, for each $k$, there are examples whose weights decrease roughly exponentially relative to the distance. Our result also leads to a characterization for perfect state transfer on unweighted unions of $\cJ(n,k)$; in particular, we obtain a simple construction that generates all unweighted unions with perfect state transfer at time $\pi/2$.

Section \ref{sec_bg} provides some background on quantum state transfer and association schemes. We then derive necessary and sufficient conditions in Section \ref{sec_wt} for a weighted graph of a Johnson scheme to exhibit perfect state transfer. Finally, \ref{sec_uw} applies the aforementioned results to unweighted unions of distance graphs.

\section{Background \label{sec_bg}}
A continuous-time quantum walk can be modeled using a real-weighted graph $X$, where each vertex represents a qubit, and each edge represents a pair of interacting qubits with coupling strength equal to the weight. In this paper, we are only concerned with quantum walks that start with an one-excitation state; under this assumption, the evolution is equivalent to the following unitary matrix acting on $\cx^{V(X)}$:
\[U(t) = \exp(-itA),\]
where $t$ denotes the time and $A$ is the weighted adjacency matrix. We say $X$ has \textsl{perfect state transfer} between $u$ and $v$ at time $\tau$ if 
\[\abs{U(\tau)_{uv}} = 1.\]

Perfect state transfer on unweighted graphs has been extensively studied. The first paper that focused on distance regular graphs appeared to be \cite{Jafarizadeh2008}, where Jafarizadeh and Sufiani paid specific attention to even cycles and hypercubes. Later on,  Coutinho, Godsil, Guo, and Vanhove \cite{Coutinho2015} thoroughly examined perfect state transfer on distance regular graphs, and more generally, graphs whose adjacency matrices lie in a Bose-Mesner algebra. Crucial to their analysis are the properties of association schemes, as we introduce now.


%



Let $I$ be the identity matrix and $J$ the all-ones matrix. An \textsl{association scheme} with $d$ classes is a set $\cA=\{\seq{A}{0}{1}{d}\}$ of $01$-matrices with the following properties:
\begin{enumerate}[(i)]
\item $A_0=I$,
\item $\sum_{r=0}^d A_r = J$,
\item $A_r^T \in \cA$ for each $r$, and 
\item $A_r A_s =A_sA_r \in \cx[\cA]$.
\end{enumerate}
We say $\cA$ is \textsl{symmetric} if $A_r=A_r^T$ for each $i$. 

The algebra $\cx[\cA]$, which has dimension $d+1$, is called the \textsl{Bose-Mesner algebra} of $\cA$. It has two bases: the Schur idempotents $\{\seq{A}{0}{1}{d}\}$, and a set of pairwise orthogonal idempotents $\{\seq{E}{0}{1}{d}\}$ with respect to the usual matrix product. Thus for each $r$, there are scalars $p_r(s)$ such that
\[A_r = \sum_{s=0}^d p_r(s) E_s.\]
These scalars $p_r(s)$ are called the \textsl{eigenvalues} of $\cA$. For more details, see Bannai and Ito \cite{Bannai1984}.

One way to construct association schemes is through distance regular graphs (see Brouwer, Cohen and Neumaier \cite{Brouwer1989}). More specifically, given a distance regular graph $X$, we define its \textsl{distance-$r$} graph $X_r$ to be the graph with the same vertex set $V(X)$ and an edge between vertices $u$ and $v$ if they are at distance $r$ in $X$. The adjacency matrices of these distance graphs form an association scheme. For example, the hypercube on $2^d$ vertices gives rise to an association scheme $\cH(d,2)$ with $d$ classes, called the \textsl{binary Hamming scheme}.


In this paper, we focus on another scheme arising from the \textsl{Johnson graph} $J(n,k,1)$. The vertex set of $J(n,k,1)$ consists of $k$-subsets of the set $\{1,2,\cdots,n\}$, and two vertices are adjacent in $J(n,k,1)$ if they differ in exactly $1$ element. The \textsl{generalized Johnson graph}, denoted $J(n,k,r)$, is defined to be distance-$r$ graph of $J(n,k,1)$. Equivalently, $J(n,k,r)$ has the same vertex set as $J(n,k,1)$, and two vertices are adjacent in $J(n,k,r)$ if they differ in exactly $r$ elements. We let $A_r$ be the adjacency matrix of $J(n,k,r)$, and call $\cJ(n,k)=\{\seq{A}{0}{1}{k}\}$ the \textsl{Johnson scheme}.

Consider a weighted graph $X$ whose adjacency matrix is Hermitian and lies in $\cx[\cJ(n,k)]$. Since $\cJ(n,k)$ is symmetric, all entries of $A(X)$ are real. For this reason, from now on, we will refer to $X$ as a \textsl{weighted graph of $\cJ(n,k)$} if $A(X)$ lives in $\re[\cJ(n,k)]$.

Unlike binary Hamming schemes, where for each positive integer $h$, there is a distance graph with perfect state transfer at time $\pi/2^h$ (see Chan \cite{Chan2013}), the only distance graphs of Johnson schemes that admit perfect state transfer are disjoint unions of edges \cite{Ahmadi2017}. Thus, it is natural to shift our attention to weighted graphs of $\cJ(n,k)$. Using the same proofs techniques in \cite[Theorem 4.1]{Godsil2011}, \cite[Theorem 4.1]{Coutinho2015} and \cite[Proposition 3.1]{Ahmadi2017}, we obtain the following necessary condition for perfect state transfer on a weighted graph of $\cJ(n,k)$.


\begin{lemma}
Let $X$ be a weighted graph of $\cJ(n,k)$. If $X$ has perfect state transfer from $u$ to $v$, then $n=2k$, and $(A_k)_{u,v}=1$.
\end{lemma}
\proof
Note that $U(t)$ is symmetric and lies in $\cx[\cJ(n,k)]$. Suppose perfect state transfer occurs between $u$ and $v$ at time $\tau$. Then $U(\tau)$ is a scalar multiple of some matrix $Q$ satisfying $Q_{uv} = 1$. Since there is exactly one class $A_r$ with $(A_r)_{uv}=1$, we must have $Q = A_r$. It follows from $A_r \circ I=0$ that $A_r$  is a permutation of order two with no fixed points. Thus, $J(n,k,r)$ is a perfect matching. On the other hand, the valency of $J(n,k,r)$ is $\binom{k}{r}\binom{n-k}{r}$, which equals $1$ if and only if $k=r$ and $n=2k$.
\qed

\section{Weighted graphs of $\cJ(2k,k)$\label{sec_wt}} 
Throughout this section, let $X$ be a weighted graph of $\cJ(2k,k)$ and write
\[A(X) = w_0A_0 + w_1A_1 + \cdots + w_mA_m,\]
for some real weights $\seq{w}{0}{1}{m}$. We characterize perfect state transfer on $X$. For $s=0,1,\cdots, k$, define 
\[\theta_s = \sum_{r=0}^m w_s p_r(s).\]
These are the eigenvalues of $X$. The following result follows from a well-known criterion for perfect state transfer on spin chains (see, for example, Kay \cite[Section 2]{Kay2010}) and that Johnson graphs are distance regular; here, we provide a direct proof using parameters of association schemes.

\begin{theorem}\label{thm_pst}
$X$ has perfect state transfer at time $\tau$ if and only if for $s=1,2,\cdots,k$, 
\[\frac{\tau(\theta_s-\theta_{s-1})}{\pi}\]
is an odd integer.
\end{theorem}
\proof
Using the spectral decomposition
\[U(\tau) = \sum_{s=0}^k e^{-i\tau \theta_{s-1}} E_s,\]
we see that $U(\tau)=\gamma A_k$ if and only if $e^{-i\tau \theta_s}= \gamma p_k(s)$ for  $s=0,1,\cdots, k$. Moreover, since $J(2k,k,1)$ is antipodal with classes of size two, by  \cite[Proposition 11.6.2]{Brouwer2012} or \cite[Lemma 4.4]{Coutinho2015}, the eigenvalues of $A_k$ are $p_k(s) = (-1)^s$. Therefore, there is some $\gamma$ such that $U(\tau)=\gamma A_k$ if and only if $\tau(\theta_s-\theta_{s-1})/\pi$ is odd for $s=1,2,\cdots,k$.
\qed

We may assume without loss of generality that perfect state transfer occurs at time $\pi$, by scaling the weights $w_j$, and that $\theta_0$ is an integer, by translating the diagonal of $A(X)$. Thus from Theorem \ref{thm_pst}, a necessary condition for $X$ to have perfect state transfer at time $\pi$ is that it has integral spectrum. To characterize such $X$, we introduce the following formula for the eigenvalues $p_r(s)$ of $\cJ(2k,k)$ \cite[Theorem 2.10]{Bannai1984}:
\begin{align*}
p_r(s) &= \sum_{j=0}^r (-1)^{r-j} \binom{k-j}{r-j}\binom{k-s}{j}\binom{k-s+j}{j}\\
&= (-1)^r \binom{k}{r} \pFq{3}{2}{-r, -k+s, k-s+1}{-k,1}{1}. 
\end{align*}
This is called the dual Hahn polynomial; for a reference, see \cite[Chapter 14]{Koekoek2010}. Setting $x=k-s$, we have
\[p_r(k-x)=\sum_{j=0}^r (-1)^{r-j}\binom{2j}{j} \binom{k-j}{r-j}\binom{x+j}{2j}.\]
Thus, the eigenvalues of $X$ are
\begin{align}\label{eqn_f(x)}
f(x)&= \sum_{r=0}^m w_r \sum_{j=0}^r p_r(k-x) \notag\\
&=\sum_{r=0}^m w_r \sum_{j=0}^r (-1)^{r-j}\binom{2j}{j} \binom{k-j}{r-j} \binom{x+j}{2j}\notag \\
&=\sum_{j=0}^m \binom{x+j}{2j} \left(\sum_{r=j}^m (-1)^{r-j} w_r \binom{2j}{j} \binom{k-j}{r-j}\right),
\end{align}
where $x=0,1,\cdots, k$. 

Let $D$ and $S$ be two $(m+1)\times (m+1)$ diagonal matrices with $D_{j,j}=\binom{2j}{j}$ and $S_{j,j} = (-1)^j$. Let $F$ be the $(k+1)\times (m+1)$ matrix with $F_{x,j} = \binom{x+j}{2j}$. Let $B$ be the $(m+1)\times (m+1)$ matrix with $B_{j,r} = \binom{k-j}{r-j}$. Let $w$ be the vector of length $m+1$ with $w_j$ as its $j$-th entry, that is, the Pascal matrix reflected about the anti-diagonal. If $c_j$ is the $j$-th entry of the vector $c=DSBSw$, then Equation (\ref{eqn_f(x)}) is equivalent to
\begin{equation}\label{eqn_f(x)withc}
f(x) = \sum_{j=0}^m c_j \binom{x+j}{2j}.
\end{equation}

Before we proceed, we remark on a connection of the above expression to integer valued polynomials. Since every matrix in $\cx[\cJ(2k,k)]$ is a polynomial in $A_1$, there is a polynomial $q(\cdot)$ such that $A(X)=q(A_1)$. As a consequence, the eigenvalues of $X$ are given by $f(x) = q(x(x+1)-k)$, and so $X$ has integral spectrum if and only if $g(y-k)$ is integer valued on $\{0\cdot 1, 1\cdot 2,\cdots,k\cdot (k+1)\}$. Using Bhargava's characterization \cite[Theorem 23]{Bhargava2000}, one can show that a polynomial is integer valued on $\{i(i+1): i\ge 0\}$ if and only if it can be written as a $\ints$-linear combination of the polynomials
\[\frac{(y-0\cdot 1)(y-1\cdot 2)\cdots (y-(j-1)j)}{(2j)!},\quad j=0,1,\cdots,\]
which, after the substitution $y=x(x+1)$, are precisely $\binom{x+j}{2j}$ for $j=0,1,\cdots$.

We now prove a characterization of weighted graphs of $\cJ(2k,k)$ with integral spectrum, using Equation (\ref{eqn_f(x)}) directly.

\begin{theorem}\label{thm_int}
$X$ has integral spectrum if and only if for $r=0,1,\cdots,m$,
\[w_r = \sum_{j=r}^m \frac{c_j}{\binom{2j}{j}} \binom{k-r}{j-r},\]
where $\seq{c}{0}{1}{m}$ are integers.
\end{theorem}
\proof
We first show that $f(x)$ in Equation (\ref{eqn_f(x)withc}) is integer valued on $\{0,1,\cdots,k\}$ if and only if $\seq{c}{0}{1}{m}$ are integers. Let $F_1$ be the matrix of the first $m+1$ rows of $F$. Since $F_1$ is unimodular, 
$F_1c$ is integral if and only if $c$ is integral. Next, from the properties of the Pascal matrix, we know that $B^{-1} = SBS$; indeed,
\begin{align*}
(BS)^2_{ji}& = \sum_{r=0}^m (-1)^{i-r} \binom{k-j}{r-j} \binom{k-r}{i-r}\\
&= \sum_{r=0}^k(-1)^{i-r} \frac{(k-j)!}{(k-i)! (i-r)!(r-j)!}\\
&= \binom{k-j}{k-i} \sum_{r=0}^k (-1)^{i-r}\binom{i-j}{i-r}\\
&=\binom{k-j}{k-i} (1-1)^{i-j}\\
&=\delta_{ij}.
\end{align*}
Therefore $w=BD^{-1}c$.
\qed

%

Recall from Theorem \ref{thm_pst} that $X$ admits perfect state transfer at time $\pi$ if and only if $f(x+1)-f(x)$ is odd, for $x=0,1,\cdots, k-1$. Using Equation (\ref{eqn_f(x)withc}), we obtain 
\[f(x+1)-f(x)= \sum_{j=1}^m c_j \binom{x+j}{2j-1}.\]

To see when the parity conditions holds, we need two useful results from number theory, due to Lucas and Kummer. For any integer $n$, let $\sigma_p(n)$ be the sum of the base-$p$ digits of $n$. The following theorem gives a formula for the largest prime power dividing a binomial coefficient; for a reference, see \cite{Dickson1966}.

\begin{theorem}\label{thm_numthy1}
Let $a$ and $b$ be non-negative integers. Let $p$ be a prime. Then
\[\nu_p \left(\binom{a}{b} \right)= \frac{\sigma_p(b)+\sigma_p (a-b) -\sigma_p(a)}{p-1}. \tag*{\sqr53}\]
\end{theorem}

\begin{theorem}\label{thm_numthy2}
Let $p$ be a prime. Let $a_n \cdots a_1a_0$ and $b_n\cdots b_1b_0$ be the representations of $a$ and $b$ in base $p$, respectively. Then
\[\binom{a}{b} \equiv \binom{a_n}{b_n}\cdots \binom{a_1}{b_1} \binom{a_0}{b_0} \pmod{p}. \tag*{\sqr53}\]
\end{theorem}

Our first result says that perfect state transfer cannot occur, regardless of the weights, if $m$ is too small.

\begin{lemma}
Let $\ell$ be the integer such that $2^{\ell} \le k <2^{\ell+1}$. If $m< 2^{\ell}$, then $X$ does not have perfect state transfer.
\end{lemma}
\proof
Let $g(x) = f(x+1)-f(x)$. We show that $g(2^{\ell}-1)$ is even. The coefficient of $c_j$ in $g(2^{\ell}-1)$ is 
\[\binom{2^{\ell} +(j-1)}{2(j-1)+1}.\]
By Theorem \ref{thm_numthy1},
\begin{align*}
\nu_2\left(\binom{2^{\ell} +(j-1)}{2(j-1)+1}\right) 
&= \sigma_2(2(j-1)+1) +\sigma_2(2^{\ell} - j) - \sigma_2(2^{\ell}+j-1)\\
&=\sigma_2(j-1)+1+\sigma_2(2^{\ell}-1-(j-1)) - 1-\sigma_2(j-1)\\
&\ge \sigma_2(2^{\ell}-j),
\end{align*}
which is at least $1$ for $j<2^{\ell}$. Therefore 
\[g(2^{\ell}-1) = \sum_{j=1}^m \binom{2^{\ell}+j-1}{2(j-1)+1}\]
is even.
\qed

On the other hand, if $m\ge 2^{\lfloor{\log_2(k)} \rfloor}$, then there exist weights $\seq{w}{0}{1}{m}$ so that $X$ has perfect state transfer.


\begin{theorem}\label{thm_par}
Let $X$ be a weighted graph of $\cJ(2k,k)$ with spectrum
\[f(x) = \sum_{j=0}^m c_j \binom{x+j}{2j},\quad x=0,1,\cdots, k-1.\]
If $2^{\ell}\le m\le k<2^{\ell+1}$, then $X$ has perfect state transfer if and only if $c_j$ is odd for $j\in \{2^0,2^1\cdots, 2^{\ell}\}$ and even otherwise.
\end{theorem}
\proof
It suffices to consider $X$ for which $\seq{c}{1}{2}{m}$ are integers. Let $G$ be the $k\times m$ matrix with $G_{x,j} = \binom{x+j}{2j-1}$, where $x=0,1,\cdots,k-1$ and $j=1,2,\cdots,m$. Then for each $x$,
\[f(x+1)-f(x) = e_x^TG\pmat{c_1\\ \vdots \\ c_m}.\]
We will work with $G$ over $\ints_2$ in the following discussion.

First, let $G_1$ be the principal $2^{\ell}\times 2^{\ell}$ submatrix of $G$. Since $G_1$ is lower triangular with $1$'s down the diagonal, its rank is $2^{\ell}$. Next, let $G_2$ be the matrix of the first $2^{\ell}$ columns of $G$. We show that the $(x,j)$-entry of $G_2$ equals the $(2(2^{\ell}-1)-x, j)$-entry of $G_2$. Indeed, Theorem \ref{thm_numthy1} gives
\[\nu_2 \left(\binom{x+j}{2j-1}\right)=\sigma_2(2j-1) + \sigma_2(x-j+1) - \sigma_2(x+j),\]
while
\begin{align*}
&\nu_2\left(\binom{2(2^{\ell}-1)-x+j}{2j-1}\right)\\
=&\sigma(2j-1) +\sigma_2(2^{\ell}-1 - (x+j) -\sigma_2 (2^{\ell+1}-1-(x-j+1))\\
=& \sigma_2(2j-1)+\ell-\sigma_2(x+j)-\ell+\sigma_2(x-j+1)\\
=&\nu_2 \left(\binom{x+j}{2j-1}\right).
\end{align*}
It follows that
\[2^{\ell}=\rk(G_1) = \rk\pmat{G_1  &\one \\ \zero  & \zero} =\rk\pmat{G_2 & | &\one }= \rk\pmat{G & | &\one }.\]
Hence, there is a unique solution to
\begin{equation}\label{eqn_c}
G\pmat{c_1\\ \vdots \\ c_m} = \one,
\end{equation}
that is,
\[\pmat{c_1\\ \vdots \\ c_{2^{\ell}}} = G_1^{-1}\one,\quad \pmat{c_{2^{\ell}+1}\\ \vdots \\ c_m} = \zero.\] 
Finally, define the new matrix
\[G_3 = \pmat{Ge_{2^0} & | &Ge_{2^1} & | &\cdots &|& Ge_{2^{\ell}}}.\]
We show that each row of $G_3$ has exactly one none-zero entry. To see this, let $a\in\{1,2,\cdots, 2^{\ell}\}$, and suppose
\[\binom{x+2^a}{2^{a+1}-1}\equiv 1 \pmod {2}.\]
By Theorem \ref{thm_numthy2}, the last $a+1$ digits in the binary representation of $x+2^a$ are all $1$. Thus $2^{a+1}$ divides $x+2^a +1$, that is, $a=\nu_2(x+1)$. Therefore, $G_3\one =\one$. Since Equation (\ref{eqn_c}) has only one solution, it holds if and only if $c_j$ is $1$ for $j\in \{2^0,2^1\cdots, 2^{\ell}\}$, and $0$ otherwise.
\qed

Theorem \ref{thm_pst} and Theorem \ref{thm_par} lead to the following characterization of perfect state transfer on weighted graphs of the Johnson scheme.

\begin{theorem}\label{thm_main}
Let $X$ be a weighted graph of $\cJ(n,k)$. Then $X$ has perfect state transfer at time $\tau$ if and only if  $n=2k$, $m\ge 2^{\lfloor{\log_2(k)} \rfloor}$, and there are integers $\seq{c}{1}{2}{m}$ such that 
\begin{enumerate}[(i)]
\item $c_j$ is odd if and only if $j$ is a power of $2$, for $j=1,2,\cdots,m$, and 
\item for  $r=1,2,\cdots,m$, 
\[w_r = \frac{\pi}{\tau} \sum_{j=r}^m \frac{c_j}{\binom{2j}{j}} \binom{k-r}{j-r}. \tag*{\sqr53}\]
\end{enumerate}
\end{theorem}

\begin{figure}[h] 
\includegraphics[width=13cm]{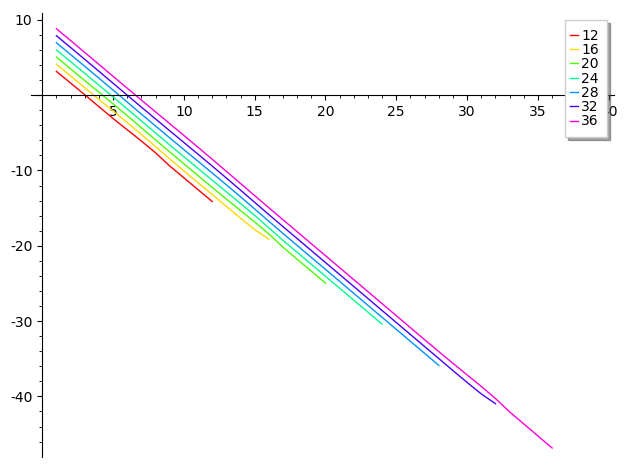}
\caption{Plot of $\log(w_j)$ for $k=12,16,\cdots,36$}
\label{fig_logplot}
\end{figure}

For example, we may set $k=m$, with $c_j=3$ if $j$ is a power of $2$ and $c_j=2$ otherwise. Figure \ref{fig_logplot} plots the corresponding weights in their logarithm, for $k=12,16,\cdots,36$.

\section{Graphs of $\cJ(2k, k)$ \label{sec_uw}}
In \cite{Ahmadi2017}, Ahmadi et al asked which unions of the distance graphs of $\cJ(2k,k)$ have perfect state transfer. We give an answer based on Theorem \ref{thm_main}.

Throughout this section, let $\cJ(2k,k) = \{\seq{A}{0}{1}{k}\}$, and let $X$ be a graph of $\cJ(2k,k)$ with adjacency matrix $\sum_{r=1}^k w_r A_r$, where $w_r\in\{0,1\}$. As before, let $f(0),\cdots,f(k)$ denote the eigenvalues of $X$. According to \cite[Theorem 4.3]{Coutinho2015}, if $X$ has perfect state transfer at time $\tau$, then $\pi/\tau$ is an odd multiple of  
\[\gcd\{f(x+1)-f(x): x=0,1,\cdots,k-1\}.\]
Thus, as early as perfect state transfer may occur, it must also occur at time $\pi/2^h$ for some non-negative integer $h$.

\begin{theorem}\label{thm_uw}
$X$ has perfect state transfer if and only if there is an integer $h\ge 1$ such that for $r=1,2,\cdots,k$,
\[w_r\equiv \sum_{j=1}^k v_j\binom{k-r}{j-r} \pmod{2^h},\]
where 
\[v_j \equiv 
\begin{cases}
2^{h-1} \pmod{2^h},\quad \sigma_2(j)=1\\
0 \pmod{2^{h-1}}, \quad \sigma_2(j)=2\\
\cdots\\
0 \pmod{2},\quad \sigma_2(j)=h.
\end{cases}\]
Moreover, if the above condition holds, then perfect state transfer occurs at time $\tau=\pi/2^h$, and possibly at time $\tau/z$ for some odd integer $z\ge 3$.
\end{theorem}
\proof
We reuse the matrices defined in Section \ref{sec_wt}, but with different indices. In the following setup, $j, r=1,2,\cdots,k$, and $x=0,1,\cdots,k-1$. Let $D$ and $S$ be two $k\times k$ diagonal matrices with $D_{j,j}=\binom{2j}{j}$ and $S_{j,j} = (-1)^j$. Let $B$ be the $k\times k$ matrix with $B_{j,r} = \binom{k-j}{r-j}$. Let $G$ be the $k\times k$ matrix with $G_{x,j} = \binom{x+j}{2j-1}$. Then
\[f(x+1)-f(x) = e_x^T GDSBSw.\]
Note that the entries of $D$ are even. By Theorem \ref{thm_main}, $X$ has perfect state transfer at time $\pi/2^h$ if and only if 
\begin{equation}\label{eqn_pst}
\left(\frac{1}{2}D\right) SBS w = 2^{h-1} c
\end{equation}
for some vector $c$ with
\[c_j \equiv  \begin{cases}
1\pmod{2},\quad \text{  $j$ is a power of $2$}\\
0\pmod{2},\quad \text{ otherwise}.
\end{cases}\]
Using the fact that 
\[D_{jj} =\binom{2j}{j} \equiv 2^{\sigma_2(j)} \pmod{2^{\sigma_2(j)+1}},\]
we can convert Equation (\ref{eqn_pst}) to
\begin{equation}\label{eqn_v}
SBSw \equiv v \pmod{2^h},
\end{equation}
where 
\[v_j \equiv 
\begin{cases}
2^{h-1} \pmod{2^h},\quad \sigma_2(j)=1\\
0 \pmod{2^{h-1}}, \quad \sigma_2(j)=2\\
\cdots\\
0 \pmod{2},\quad \sigma_2(j)=h.
\end{cases}\]
Finally, since $(SB)^2=I$, Equation (\ref{eqn_v}) holds if and only if $w\equiv Bv \pmod{2^h}$.

Once the weights $w_r$'s are determined, one may compute $c$ and obtain the greatest common divisor of its entries, say $z$. Then perfect state transfer occurs at time $\pi/(2^hz)$, as well as $\pi/2^h$.
\qed

We immediately obtain the following result, which constructs all graphs of the Johnson scheme with perfect state transfer at time $\pi/2$.

\begin{theorem}\label{thm_uwpi/2}
$X$ has perfect state transfer at time $\pi/2$ if and only if for some subset $T$ satisfying
\[\{2^0, 2^1,\cdots, 2^{\lfloor\log_2(k)\rfloor}\}\sbs T \sbs\{1,2,\cdots,k\},\]
we have
\[w_r \equiv \sum_{j\in T} \binom{k-r}{j-r} \pmod{2},\]
where $r=1,2,\cdots,k$.
\qed
\end{theorem}

\section{Future work}
While Theorem \ref{thm_uw} gives a characterization for graphs of $\cJ(2k,k)$ with perfect state transfer, it is not easy to generate all of them as Corollary \ref{thm_uwpi/2} does for $h=1$. However, there exist examples with perfect state transfer at time $\pi/4$; the smallest one lives in $\cx[\cJ(62,31)]$, and is a union of the distance-$r$ graphs where $r\in\{7, 11, 13, 14, 15, 19, 21, 22, 23, 25, 26, 27, 28, 29, 30\}$. Thus, a constructive characterization based on Theorem \ref{thm_uw} is desirable, and it may help to look at the number-theoretic properties of the Pascal matrix.

Questions (i) and (ii) in Section \ref{sec_intro} are also interesting directions to work on. In particular, some of our techniques can be applied to other association schemes, such as the Hamming scheme. It is known that there are distance graphs \cite{Chan2013}, as well as weighted unions of $X_1$ and $X_2$ \cite{Christandl2017}, of $\cH(d,2)$ that have perfect state transfer, but a complete list is yet to be determined. We plan on following up on these questions.

\bibliographystyle{amsplain}
\bibliography{../../../../Bibtex/cqw.bib}

\end{document}